УДК 681.51 DOI: 10.17587/mau.22….

**А. А. Бобцов[1],** д-р техн. наук, проф., bobtsov@mail.ru,
**Н. А. Николаев[1],** канд. техн. наук, доц., nikona@yandex.ru,
**Р. Ортега[2],** д-р техн. наук, проф., romeo.ortega@itam.mx,
**О. В. Слита[1],** канд. техн. наук, доц., o-slita@yandex.ru,
**О. А. Козачёк[1],** аспирант, oakozachek@mail.ru
[1] *Университет ИТМО, г. Санкт-Петербург,*
[2] *Instituto Tecnológico Autónomo de México, México*

# Адаптивный наблюдатель переменных состояния линейной нестационарной системы с частично неизвестными параметрами матрицы состояния и вектора входа[1]

*В статье рассматривается проблема синтеза адаптивного наблюдателя переменных состояния линейной нестационарной одноканальной динамической системы. Предполагается, что сигнал управления и выходная переменная измеряемы. При этом допускается, что матрица состояния объекта управления содержит известные переменные и неизвестные постоянные параметры, а матрица (вектор) управления неизвестна. Синтез наблюдателя основан на методе GPEBO (обобщенный наблюдатель, основанный на оценке параметров), предложенном в работе [1]. Синтез адаптивного наблюдателя предусматривает предварительную параметризацию исходной системы и преобразование ее к линейной регрессионной модели с дальнейшей идентификацией неизвестных параметров. Для идентификации неизвестных постоянных параметров был использован классический алгоритм оценки – метод наименьших квадратов с фактором забывания (forgetting factor). Данный подход хорошо себя зарекомендовал в случаях, когда известный регрессор является «частотно бедным» (то есть спектральный состав регрессора содержит менее r/2 гармоник, где r – число неизвестных параметров) или не удовлетворяет так называемому условию незатухающего возбуждения. Для иллюстрации работоспособности предложенного метода в статье представлен пример. Был рассмотрен нестационарный объект второго порядка с четырьмя неизвестными параметрами. Была произведена параметризация исходной динамической модели и получена линейная статическая регрессия, содержащая шесть неизвестных*

---

*параметров (включая вектор неизвестных начальных условий переменных состояния системы). Был синтезирован адаптивный наблюдатель и представлены результаты компьютерного моделирования, иллюстрирующие достижение заданной цели. Основным отличием от результатов, опубликованных ранее в работе [2], является новое допущение о том, что линейная нестационарная система содержит не только неизвестные параметры в матрице состояния, но и матрица (вектор) по управлению содержит неизвестные постоянные коэффициенты.*

***Ключевые слова:*** *адаптивный наблюдатель, линейная нестационарная система, линейная регрессионная модель*

## Введение

В задачах управления сложными динамическими системами получение информации о состоянии системы осуществляется с помощью первичных измерительных преобразователей (датчиков). При этом не всегда удается разместить набор измерительных средств, который позволит измерить весь вектор состояния объекта, поэтому для оценки недоступных прямым измерениям переменных состояния используются наблюдатели. Для случаев, когда объект описывается линейной динамической системой с постоянными параметрами, в настоящее время разработаны эффективные методы синтеза наблюдателей. При этом стоит отметить, что не всегда возможно достаточно точно описать поведение объекта с помощью моделей с постоянными параметрами в связи с тем, что параметры могут изменяться с течением времени. Это обусловлено различными как внутренними, так и внешними факторами, такими, например, как старение и соответствующее ухудшение параметров элементов системы, воздействие температуры, изменение массо-габаритных параметров в процессе функционирования. С примерами использования линейных нестационарных систем к задачам механики можно ознакомиться, например в [3]. Таким образом, более точно поведение сложной динамической системы может быть описано с помощью дифференциальных уравнений с нестационарными параметрами, поэтому в настоящее время большую актуальность имеют исследования, направленные на синтез наблюдателей для нестационарных систем.

Проблема синтеза наблюдателей для линейных нестационарных систем не является новой, однако интерес исследователей к данной проблеме не угасает. Так, например, в [4] рассматривается задача синтеза оптимальных эллипсоидных наблюдателей и алгоритмов идентификации, обеспечивающих наилучшие эллипсоидные оценки состояния системы и неизвестных параметров. Задачи построения наблюдателей решаются различными подходами, одним из которых является сведение исходной математической модели к линейной регрессионной (см., например, [2], [5]) с последующей идентификацией ее параметров.

Проблема оценки вектора состояния линейной нестационарной системы является актуальной как при синтезе управления, так и имеет самостоятельное значение. Например, информация, получаемая с наблюдателя, может использоваться для построения устройств контроля технического состояния объекта [6]…[8].

В задачах оценивания вектора состояния линейных нестационарных систем с частично неизвестными параметрами обычно рассматриваются динамические системы с различными допущениями относительно матриц с нестационарными параметрами вида

$$\dot{\mathbf{x}}(t) = \mathbf{A}(t)\mathbf{x}(t) + \mathbf{B}(t)u(t) + \mathbf{W}(t)\boldsymbol{\theta}, \qquad (1)$$
$$y(t) = \mathbf{C}^{\mathrm{T}}(t)\mathbf{x}(t),$$

где $\mathbf{x}(t)$, $u(t)$, $y(t)$ векторы состояния, управления и выхода соответствующей размерности (при этом в различных постановках рассматриваются как системы с одним входом и одним выходом (SISO), как, например в [9], [10] так и многомерные системы (MIMO), как, например, в [11]), $\mathbf{A}(t)$, $\mathbf{B}(t)$, $\mathbf{C}^T(t)$ – известные матрицы с нестационарными коэффициентами соответствующей размерности, $\boldsymbol{\theta}$ – вектор неизвестных параметров и $\mathbf{W}(t)$ – матрица известных сигналов соответствующей размерности.

В качестве цели ставится задача оценки вектора состояния и вектора неизвестных параметров по измерениям $u(t)$, $y(t)$ и $\mathbf{W}(t)$.

При отсутствии вектора неизвестных параметров $\boldsymbol{\theta}$ задача сводится к классической проблеме оценки вектора состояния, например [12], для которой алгоритм оценки ищется в виде

$$\dot{\hat{\mathbf{x}}}(t) = \mathbf{F}(t)\hat{\mathbf{x}}(t) + \mathbf{G}(t)y(t) + \mathbf{H}(t)u(t),$$

где $\mathbf{F}(t), \mathbf{G}(t), \mathbf{H}(t)$ – соответствующие матрицы.

В настоящее время проблемам синтеза наблюдателей для неизвестного вектора состояния нестационарной системы посвящено достаточно большое число работ, при этом вводятся различные предположения и допущения относительно матриц описания объекта. Так, например в [13] задача синтеза наблюдателя решается в предположении, что матрица объекта $\mathbf{A}(t)$ представлена в канонической форме. В данной работе предлагается развитие подхода, предложенного ранее в работе [2], где рассматривалась линейная нестационарная система с частично неизвестными параметрами, в отличие от полученного в ней результата в данной работе выполнено расширение на случай, когда матрица (вектор) управления содержит неизвестные постоянные параметры.

## Постановка задачи

Рассматривается линейная нестационарная система с одним входом и одним выходом (SISO) вида:

$$\dot{\mathbf{x}}(t) = \mathbf{A}(t)\mathbf{x}(t) + \mathbf{k}\mathbf{C}^{\mathrm{T}}(t)\mathbf{x}(t) + \mathbf{b}u(t), \mathbf{x}(0) = \mathbf{x}_0 \in \mathbb{R}^n, t \geq 0,$$
$$y(t) = \mathbf{C}^{\mathrm{T}}(t)\mathbf{x}(t), \qquad (2)$$

где $\mathbf{x}(t) \in \mathbb{R}^n$ – неизмеряемый вектор состояния, $u(t) \in \mathbb{R}$ – известный входной сигнал, $y(t)$ – измеряемый выходной сигнал, матрицы $\mathbf{A}(t), \mathbf{C}^{\mathrm{T}}(t)$ являются известными матрицами с нестационарными параметрами, $\mathbf{k}$ и $\mathbf{b}$ - векторы соответствующей размерности. Предполагается, что $\mathbf{A}(t)$ и $\mathbf{C}(t)$ известны и ограничены, а векторы $\mathbf{k} \in \mathbb{R}^n$ и $\mathbf{b} \in \mathbb{R}^n$ постоянны и неизвестны.

На основе этих данных в статье разработан адаптивный наблюдатель:

$$\dot{\boldsymbol{\chi}}(t) = \mathbf{F}\big(\boldsymbol{\chi}(t), u(t), y(t)\big),$$

$$\begin{bmatrix} \hat{\mathbf{x}}(t) \\ \hat{\mathbf{k}}(t) \\ \hat{\mathbf{b}}(t) \end{bmatrix} = \mathbf{S}(\chi(t), u(t), y(t)),$$

где $\chi(t) \in \mathbb{R}^{n_\chi}$ такой, что все сигналы ограничены, а также обеспечивается сходимость оценок переменных состояния и неизвестных постоянных параметров к их реальным значениям:

$$\hat{\mathbf{x}}(t) = \mathbf{x}(t), \hat{\mathbf{k}}(t) = \mathbf{k}, \hat{\mathbf{b}}(t) = \mathbf{b},$$

для всех $\mathbf{x}_0 \in \mathbb{R}^n$, $\chi(t) \in \mathbb{R}^{n_\chi}$.

В соответствии со стандартным подходом в теории наблюдателей, траектории входа и состояния полагаются ограниченными.

Рассмотрим типовые предположения, которые накладываются на матрицы исходной системы (см., например, [14]…[16]).

Предположение 1. Пара матриц $\mathbf{A}(t)$ и $\mathbf{C}^T(t)$ является обнаруживаемой, то есть существует матрица обратной связи $L(t)$ такая, что система

$$\dot{\mathbf{x}}(t) = [\mathbf{A}(t) - \mathbf{L}(t)\mathbf{C}^T(t)]\mathbf{x}(t),$$

является асимптотически устойчивой.

Предположение 2. Фундаментальная матрица, автономной системы $\dot{\mathbf{x}}(t) = \mathbf{A}_0(t)\mathbf{x}(t)$, где $\mathbf{A}_0(t) = \mathbf{A}(t) - \mathbf{L}(t)\mathbf{C}^T(t)$, удовлетворяет условию

$$\|\mathbf{\Phi}_{\mathbf{A}_0}(t, \tau)\| \leq c_1, \forall\, t \geq \tau \geq 0,$$

что означает равномерную устойчивость автономной системы (uniformly stable), см. Теорему 6.4 [16].

Предположение 3. Выполняется условие вида

$$\int_\tau^t \|\mathbf{\Phi}_{\mathbf{A}_0}(t, s)\mathbf{B}(s)\| ds \leq c_2, \forall\, t \geq \tau \geq 0,$$

являющееся необходимым и достаточным условием устойчивости типа ограниченный вход ограниченное состояние (BIBS – bounded-input-bounded-state) системы (1) с добавлением составляющей $-\mathbf{L}(t)\mathbf{C}^T(t)$ и при $\mathbf{k} = 0$.

## Основной результат

В данной работе на первом шаге решается задача параметризации, то есть сведение задачи оценки параметров для исходной динамической модели к идентификации параметров для линейной статической регрессионной модели. Вторым шагом выполняется оценка неизвестных постоянных параметров линейной регрессионной модели, которая, как известно, может решаться различными методами в зависимости от того, какие условия возбуждения накладываются на регрессор, см., например, [5], [17], [18].

**Теорема 1.** Рассмотрим динамическую систему вида

$$\dot{\boldsymbol{\xi}}(t) = \mathbf{A}_0(t)\boldsymbol{\xi}(t) + \mathbf{L}(t)y(t), \ \boldsymbol{\xi}(0) = \mathbf{0}_{n\times 1}, \quad (3)$$

$$\dot{\boldsymbol{\eta}}(t) = \boldsymbol{A}_0(t)\boldsymbol{\eta}(t) + \mathbf{I}y(t), \quad \boldsymbol{\eta}(0) = \mathbf{0}_{n\times n}, \quad (4)$$

$$\dot{\boldsymbol{\zeta}}(t) = \mathbf{A}_0(t)\boldsymbol{\zeta}(t) + \mathbf{I}u(t), \quad \boldsymbol{\zeta}(0) = \mathbf{0}_{n\times n}, \quad (5)$$

$$\dot{\boldsymbol{\Phi}}(t) = \mathbf{A}_0(t)\boldsymbol{\Phi}(t), \quad \boldsymbol{\Phi}(0) = \mathbf{I}_{n\times n},$$

где $\mathbf{A}_0(t)$ и $\mathbf{L}(t)$ удовлетворяют допущениям 1…3.

Тогда исходная динамическая система (2) может быть преобразована к линейной регрессионной модели вида

$$z(t) = \boldsymbol{\Psi}(t)\boldsymbol{\Theta}, \quad (6)$$

где функция $z(t) = y(t) + \mathbf{C}^\mathrm{T}(t)\boldsymbol{\xi}(t)$ является известной, вектор известных функций $\boldsymbol{\Psi}(t) = [\mathbf{C}^\mathrm{T}(t)\boldsymbol{\Phi}(t) \quad \mathbf{C}^\mathrm{T}(t)\boldsymbol{\eta}(t) \quad \mathbf{C}^\mathrm{T}(t)\boldsymbol{\zeta}(t)]$ и $\boldsymbol{\Theta} = [\boldsymbol{\theta} \quad \mathbf{k} \quad \mathbf{b}]^\mathrm{T}$ – вектор неизвестных параметров.

*Доказательство.*

Рассмотрим уравнение ошибки вида

$$\mathbf{e}(t) = \boldsymbol{\xi}(t) + \boldsymbol{\eta}(t)\mathbf{k} + \boldsymbol{\zeta}(t)\mathbf{b} - \mathbf{x}(t), \quad (7)$$

тогда для производной от ошибки имеем

$$\dot{\mathbf{e}}(t) = \dot{\boldsymbol{\xi}}(t) + \dot{\boldsymbol{\eta}}(t)\mathbf{k} + \dot{\boldsymbol{\zeta}}(t)\mathbf{b} - \dot{\mathbf{x}}(t) = \mathbf{A}_0(t)\boldsymbol{\xi}(t) + \mathbf{L}(t)y(t) +$$
$$+\mathbf{A}_0(t)\boldsymbol{\eta}(t)\mathbf{k} + \mathbf{I}y(t)\mathbf{k} + \mathbf{A}_0(t)\boldsymbol{\zeta}(t)\mathbf{b} + \mathbf{I}u(t)\mathbf{b} -$$
$$- \big(\mathbf{A}_0(t) + \mathbf{L}(t)\mathbf{C}^T(t)\big)\mathbf{x}(t) - \mathbf{k}\mathbf{C}^\mathrm{T}(t)\mathbf{x}(t) - \mathbf{b}u(t) =$$
$$= \mathbf{A}_0(t)\boldsymbol{\xi}(t) + \mathbf{A}_0(t)\boldsymbol{\eta}(t)\mathbf{k} + \mathbf{A}_0(t)\boldsymbol{\zeta}(t)\mathbf{b} - \mathbf{A}_0(t)\mathbf{x}(t) = \mathbf{A}_0\mathbf{e}(t)$$
$$\dot{\mathbf{e}}(t) = \mathbf{A}_0\mathbf{e}(t).$$

Решение дифференциального уравнения $\dot{\mathbf{e}}(t) = \mathbf{A}_0\mathbf{e}(t)$ имеет вид

$$\mathbf{e}(t) = \boldsymbol{\Phi}(t)\boldsymbol{\theta}, \quad (8)$$

где $\boldsymbol{\theta} = \mathbf{e}(0)$ и в случае, если начальные условия для динамической системы (3)…(5) выбраны нулевыми, то $\mathbf{e}(0) = -\mathbf{x}(0)$.

Далее, после подстановки (8) в (7) имеем

$$\mathbf{x}(t) - \boldsymbol{\xi}(t) = \boldsymbol{\eta}(t)\mathbf{k} + \boldsymbol{\zeta}(t)\mathbf{b} - \boldsymbol{\Phi}(t)\boldsymbol{\theta}. \qquad (9)$$

Домножив слева левую и правую части выражения (9) на $\mathbf{C}^{\mathrm{T}}(t)$ получаем линейную регрессионную модель вида (6).

Если выполняется условие неисчезающего возбуждения

$$\alpha_2 \mathbf{I} \leq \int_{t_0}^{t_0+\delta} \boldsymbol{\Psi}(\tau)\boldsymbol{\Psi}^T(\tau) d\tau \leq \alpha_1 \mathbf{I}, \text{для всех } t_0 > 0,$$

где $\alpha_1$, $\alpha_2$ и $\delta$ – положительные константы, то для нахождения неизвестных параметров линейной регрессионной модели можно использовать различные подходы, в том числе градиентный алгоритм идентификации [19], метод динамического расширения регрессора и смешивания [17], [18] и др.

В данной работе для оценки параметров линейной регрессионной модели предлагается использовать классический алгоритм оценки – метод наименьших квадратов с фактором забывания (forgetting factor) [20], [21] вида

$$\dot{\widehat{\boldsymbol{\theta}}} = \gamma \mathbf{F}(t)\boldsymbol{\Psi}^{\mathrm{T}}(t)\big(z(t) - \boldsymbol{\Psi}(t)\widehat{\boldsymbol{\theta}}\big),$$

$$\dot{\mathbf{F}} = \begin{cases} -\gamma \mathbf{F}(t)\boldsymbol{\Psi}^{\mathrm{T}}(t)\boldsymbol{\Psi}(t)\mathbf{F}(t) + \beta \mathbf{F}(t), & \text{если } \|\mathbf{F}(t)\| \leq \mathrm{M} \\ \mathbf{0} & \text{если } \|\mathbf{F}(t)\| > \mathrm{M}' \end{cases}$$

где $\mathbf{F}(0) = \frac{1}{f_0}\mathbf{I}$, где $\mathbf{I}$ – единичная матрица соответствующей размерности, настраиваемые параметры $\gamma > 0$, $\beta > 0$, $f_0 \geq 0$, $\mathrm{M} > 0$.

Таким образом, оценку вектора состояния динамической системы (2) можно получить из (9) в виде

$$\widehat{\mathbf{x}}(t) = \boldsymbol{\xi}(t) - \boldsymbol{\Phi}(t)\begin{bmatrix}\Theta_1\\\Theta_2\end{bmatrix} + \boldsymbol{\eta}(t)\begin{bmatrix}\Theta_3\\\Theta_4\end{bmatrix} + \boldsymbol{\zeta}(t)\begin{bmatrix}\Theta_5\\\Theta_6\end{bmatrix}$$

**Результаты моделирования**

Моделирование проводилось при следующих параметрах системы (2)

$$\mathbf{A}(t) = \begin{bmatrix} 1{,}8 + \sin(0{,}5t) & -1 \\ 5.2 + \cos(2t) + 0{,}5\sin(t) & -4 \end{bmatrix}, \mathbf{b} = \begin{bmatrix}1\\2\end{bmatrix},$$

$$\boldsymbol{k} = \begin{bmatrix} -1 \\ -3 \end{bmatrix}, \mathbf{C}(t) = \begin{bmatrix} 1 \\ 0 \end{bmatrix}.$$

Используя вектор $\mathbf{L}(t) = \begin{bmatrix} 0{,}8 + 0{,}5\sin(0{,}5t) \\ 0{,}2 + \cos(2t) \end{bmatrix}$, имеем

$$\mathbf{A}_0(t) = \begin{bmatrix} 1 & -1 \\ 5 + 0{,}5\sin(t) & -4 \end{bmatrix}.$$

Моделирование проводилось при начальных условиях $\mathbf{x}(0) = \begin{bmatrix} 3 \\ -2 \end{bmatrix}$. Параметры алгоритма оценки были выбраны следующим образом $\alpha = 1000$, $M = 10^{12}$, $\beta = 1$, $f_0 = 0{,}1$, $u(t) = \sin(t)$.

На рисунке 1 приведены переходные процессы оценки неизвестных параметров $\widehat{\Theta} = [\hat{\mathbf{e}}(\mathbf{0}) \quad \hat{\mathbf{k}} \quad \hat{\mathbf{b}}]^{\mathrm{T}}$. На рисунках 2…7 приведены результаты моделирования ошибок оценивания неизвестных параметров $\widehat{\Theta}_i - \Theta_i$, где $i = \overline{1,6}$. На рисунках 8 и 9 приведены переходные процессы по ошибке оценки переменных состояния нестационарной системы $\hat{x}_i - x_i$, где $i = \overline{1,2}$.

Результаты моделирования, приведенные на рисунках 1…9, иллюстрируют работоспособность предложенного в работе алгоритма идентификации переменных состояния линейной нестационарной системы (2) с частично неизвестными параметрами и неизвестной постоянной матрицей управления.

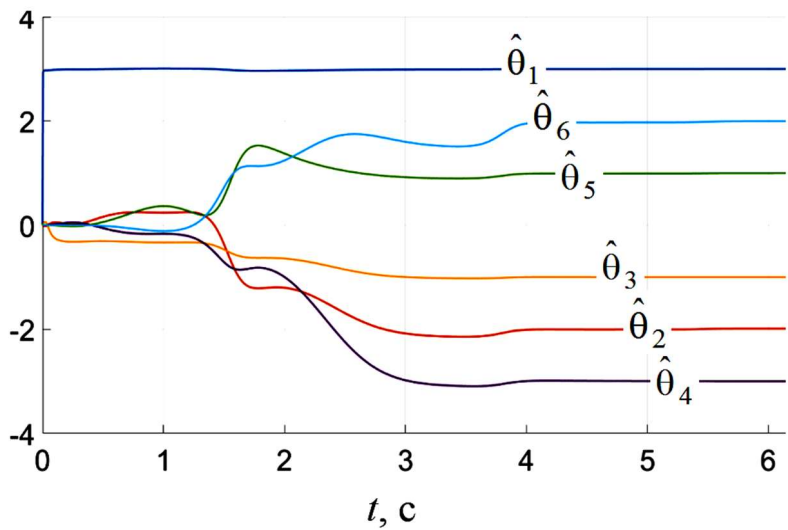

**Рис. 1.** Переходные процессы по оценки неизвестных параметров
**Fig. 2.** Transient processes of estimation of unknown parameters

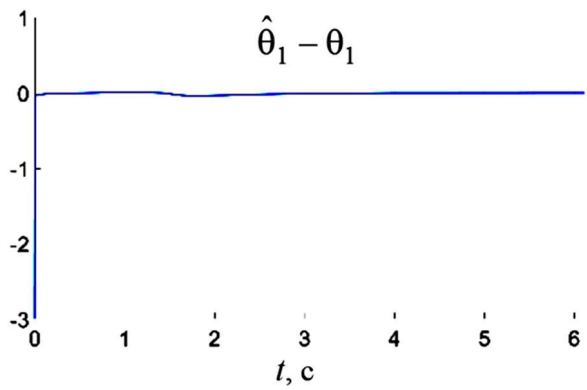

Рис. 2. Переходный процесс
по ошибке $\hat{\theta}_1 - \theta_1$

**Fig. 2.** Transient process of error $\hat{\theta}_2 - \theta_2$

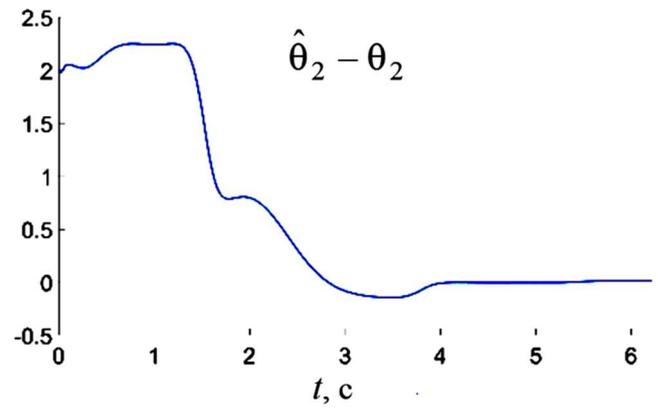

Рис. 3. Переходный процесс
по ошибке $\hat{\theta}_2 - \theta_2$

**Fig. 3.** Transient process of error $\hat{\theta}_2 - \theta_2$

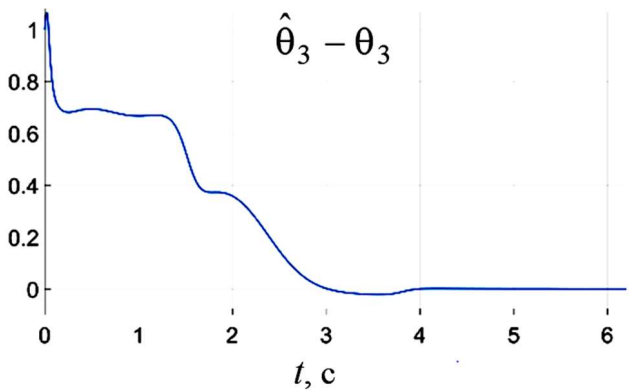

Рис. 4. Переходные процесс
по ошибке $\hat{\theta}_3 - \theta_3$

**Fig. 4.** Transient process of error $\hat{\theta}_3 - \theta_3$

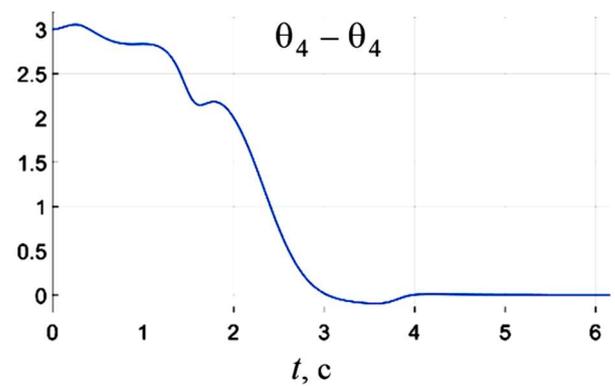

Рис. 5. Переходный процесс
по ошибке $\hat{\theta}_4 - \theta_4$

**Fig. 5.** Transient process of error $\hat{\theta}_4 - \theta_4$

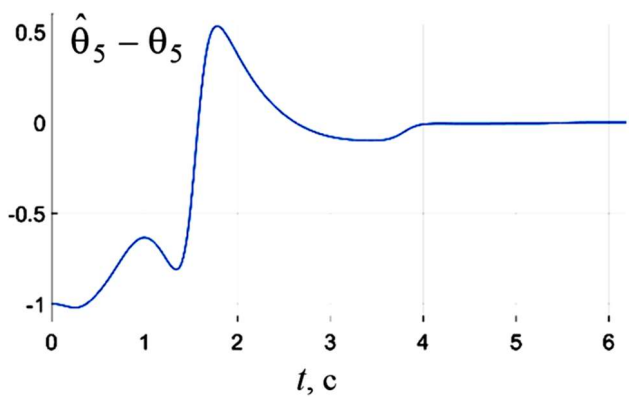

Рис. 6. Переходной процесс
по ошибке $\hat{\theta}_5 - \theta_5$

**Fig. 6.** Transient process of error $\hat{\theta}_5 - \theta_5$

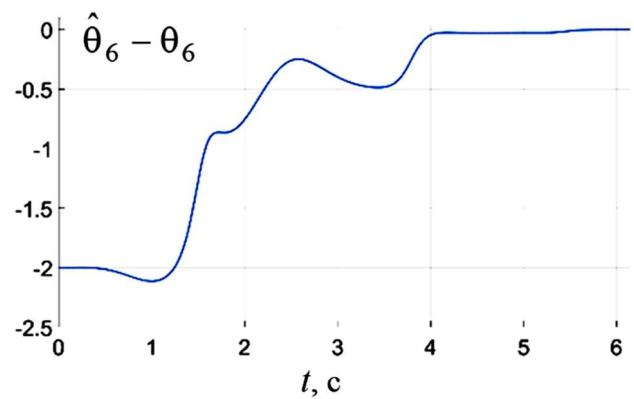

Рис. 7. Переходной процесс
по ошибке $\hat{\theta}_6 - \theta_6$

**Fig. 7.** Transient process of error $\hat{\theta}_6 - \theta_6$

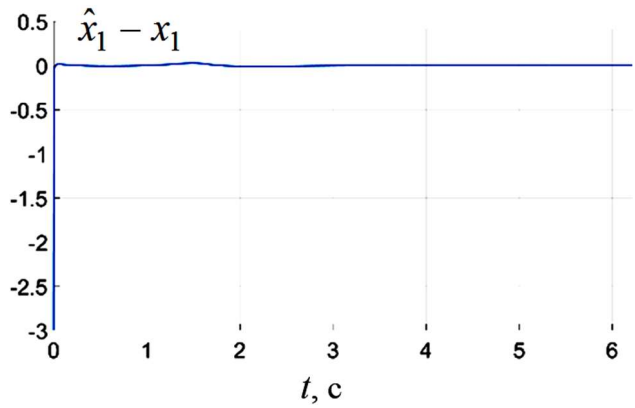 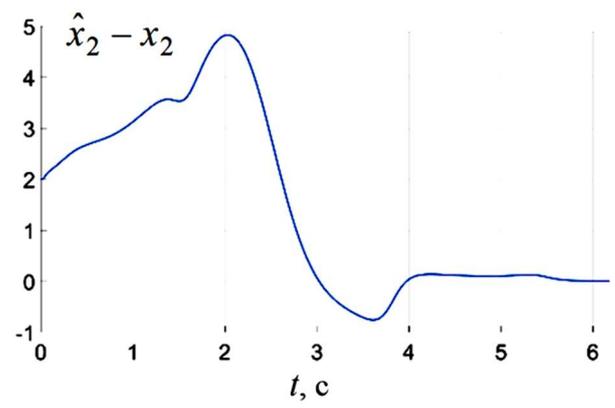

**Рис. 8.** Переходной процесс по ошибке $\hat{x}_1 - x_1$

**Fig. 8.** Transient process of error $\hat{x}_1 - x_1$

**Рис. 9.** Переходной процесс по ошибке $\hat{x}_2 - x_2$

**Fig. 9.** Transient process of error $\hat{x}_2 - x_2$

## Заключение

В статье предложен адаптивный наблюдатель состояния линейной нестационарной системы с частично неизвестными постоянными параметрами. Задача решена в предположении, что измеряется только выходная переменная, и матрицы состояния и управления содержат неизвестные постоянные значения. Результаты математического моделирования иллюстрируют работоспособность предложенного алгоритма.

# Список литературы

# Adaptive state observer for linear time-variant system with partially unknown state matrix and input matrix parameters

**A. A. Bobtsov[1]**, *bobtsov@mail.ru*,

**N. A. Nikolaev[1]**, *nikona@yandex.ru*,

**R. Ortega[2]**, *romeo.ortega@itam.mx*,

**O. V. Slita[1]**, *o-slita@yandex.ru*,

**O. A. Kozachek[1]**, *oakozachek@mail.ru*

[1] *ITMO University, Saint-Petersburg, Russia*

[2] *Instituto Tecnológico Autónomo de México, México*

*Corresponding author*:

**Bobtsov A. A.**, Dr. Sci., Professor, ITMO University, Saint Petersburg, 197101, Russian Federation

e-mail: bobtsov@mail.ru

### Abstract

*In this paper the problem of adaptive state observer synthesis for linear time-variant SISO (single-input-single-output) dynamical system with partially unknown parameters was considered. It is assumed that the input signal and output variable of the system are measurable. It is also assumed that the state matrix of the plant contains known variables and unknown constants when the input matrix (vector) is unknown. Observer synthesis is based on GPEBO (generalized parameter estimation based observer) method proposed in [1]. Observer synthesis provides preliminary parametrization of the initial system and its conversion to a linear regression model with further unknown parameters identification. For identification of the unknown constant parameters classical estimation algorithm – least squares method with forgetting factor – was used. This approach works well in cases, when the known regressor is «frequency poor» (i.e. the regressor spectrum contains r/2 harmonics, where r is a value of the unknown parameters) or does not meet PE (persistent excitation) condition. To illustrate performance of the proposed method, an example is provided in this paper. A time-variant second-order plant with four unknown parameters was considered. Parametrization of the initial dynamical model was made. A linear static regression with six unknown parameters (including unknown state initial conditions vector) was*

*obtained. An adaptive observer was synthesized and the simulation results were provided to illustrate the purpose reached. The main difference with the results, that were published earlier in [2], is the new assumption that not only does the state matrix of the linear time-variant system contain unknown parameters, but input matrix (vector) contains unknown constant coefficients.*

*Keywords*: adaptive observer, linear time-variant system, linear regression model

**Acknowledgements:** This work was supported by Russian Science Foundation, project no. 22-21-00499, https://rscf.ru/project/22-21-00499/.

DOI: 10.17587/mau.22……………